\documentclass{article}
\usepackage{verbatim}
\usepackage{amsmath} 
\usepackage{amsthm}
\usepackage{amssymb}
\usepackage{enumerate}
\usepackage{amssymb}
\newtheorem{theorem}{Theorem}

\newtheorem{lemma}[theorem]{Lemma}

\newtheorem{claim}{Claim}[theorem]

\newtheorem{obs}[theorem]{Observation}

\theoremstyle{definition}

\newtheorem{definition}[theorem]{Definition}

\theoremstyle{remark}

\def\myheads#1;#2;{
\pagestyle{myheadings}
\markboth{{\sc\hfill #1\hfill\protect\makebox[0cm][r]{\rm\today}}}
{{\sc\protect\makebox[0cm][l]{\rm\today}\hfill #2\hfill}}
}

\newcommand{\acal}{{\mathcal A}}
\newcommand{\bcal}{{\mathcal B}}
\newcommand{\ccal}{{\mathcal C}}

\newcommand{\pcal}{{\mathcal P}}

\newcommand{\scal}{{\mathcal S}}

\newcommand{\xcal}{{\mathcal X}}

\newcommand{\setm}{\setminus}
\newcommand{\empt}{\emptyset}
\newcommand{\subs}{\subset}

\newcommand{\oo}{{{\omega}_1}}
\newcommand{\dom}{\operatorname{dom}}

\def\<{\left\langle}
\def\>{\right\rangle}
\def\cf{\operatorname{cf}}

\def\OO{{\omega}}
\def\oo{\omega_1}
\def\br#1;#2;{\bigl[ {#1} \bigr]^ {#2} }
\def\bc#1;#2;{\bigl( {#1} \bigr)^ {#2} }

\def\ooseq#1;#2;{\< {#1}_{#2}:{#2}<\oo\>}
\def\ooset#1;#2;{\{ {#1}_{#2}:{#2}<\oo\}}
\def\seq#1;#2;#3;{\< {#1}_{#2}:{#2}<#3\>}
\def\set#1;#2;#3;{\{ {#1}_{#2}:{#2}<#3\}}
\def\oseq#1;#2;{\< {#1}_{#2}:{#2}<\OO\>}
\def\oset#1;#2;{\{ {#1}_{#2}:{#2}<\OO\}}
\def\oosequ#1;#2;{\< {#1}^{#2}:{#2}<\oo\>}
\def\oosetu#1;#2;{\{ {#1}^{#2}:{#2}<\oo\}}
\def\sequ#1;#2;#3;{\< {#1}^{#2}:{#2}<#3\>}
\def\setu#1;#2;#3;{\{ {#1}^{#2}:{#2}<#3\}}
\def\osequ#1;#2;{\< {#1}^{#2}:{#2}<\OO\>}
\def\osetu#1;#2;{\{ {#1}^{#2}:{#2}<\OO\}}

\def\to{\longrightarrow}
\def\rank{\operatorname{rank}}

\def\fin#1;{\br #1;<{\omega};}

\newcommand{\restr}%
{\mathop{\hspace{0.01ex}|\hspace*{-0.02ex}{\grave{}}\hspace{0.4ex}}}

\newcommand{\rest}{\restr}

\begin{document}
\title
{Boolean algebras with prescribed topological densities\thanks{The
preparation of this paper was supported by the  
Hungarian National Foundation for Scientific Research grant no. 25745 and
JSPS grant n.o. P 98259 }}
\author{L. Soukup}
\maketitle 
\myheads{\date Boolean algebras with prescribed topological density};
{BA's with prescribed  density};
\newcommand{\ise}{\triangleleft}
\newcommand{\len}{\ell}
\newcommand{\conc}{{}^\frown\!}
\newcommand{\calb}{{\mathcal B}}
\newcommand{\calbt}{{\mathcal B}_{\mu}}
\newcommand{\calf}{{\mathcal F}}
\newcommand{\base}{\operatorname{\mathbf b}}
\newcommand{\den}{\operatorname{d}}
\newcommand{\trivalg}{{\mathbf 2}}
 
\def\ta#1;{\ell_{t}(#1)} 
\def\he#1;{\ell_{h}(#1)}
\def\tas#1;{{#1}^{\uparrow}} 
\def\hes#1;{{#1}^{\downarrow}}
\def\nt{sequence with nice tail\ } 
\def\tnt{T_{\rm nt}}

\begin{abstract}
We give a simplified proof of a theorem of M. Rabus and S. Shelah
claiming that for each cardinal ${\mu}$
there is a c.c.c Boolean algebra $\bcal$ with topological density
${\mu}$.
\end{abstract}

The {\em topological density $\den(B)$} of a Boolean algebra $B$
is the minimal cardinal ${\mu}$ such that $B\setm \{0_B\}$
can be covered by ${\mu}$ ultrafilters.

\begin{theorem}[M. Rabus and S. Shelah]\label{tm:rs}
For each cardinal ${\mu}$
there is a c.c.c Boolean algebra $\bcal$ with topological density
${\mu}$. 
\end{theorem}

The simplest way to guarantee  $\den(B)\ge {\mu}$ for some Boolean
algebra $B$ is to find
a family $\acal\subs B\setm \{0_B\}$ of size ${\mu}$
such that $a\land a'=0$ for each $\{a,a'\}\in \br \acal;2;$.
Obviously if $B$ satisfies c.c.c then this argument can not work for 
${\mu}\ge \oo$.
However, 
 instead of finding  ${\mu}$ elements of $B$ with pairwise  
empty intersections it is enough to require
that $B$ contains sequences with 
``{
prescribedly small  pairwise intersections''}: 
\begin{obs}
Let $B$ be a Boolean algebra, ${\mu}$ be a cardinal.
Assume that there is a set $X\subs B\setm \{0\}$ such that 
$X$ is not the union of finitely many centered sets and 
for each ${\nu}<{\mu}$
we have 
\begin{multline}\tag{$\dag$}\label{obs}
\forall \<x_{\xi}:{\xi}<{\nu}^+\>\subs X\ 
\forall F:\br {\nu}^+;2;\stackrel{1-1}{\longrightarrow}{\nu}^+\\
\exists \<y_{\xi}:{\xi}<{\nu}^+\>\subs X\ 
\forall \{{\xi},{\zeta}\}\in \br {\nu}^+;2;
\  y_{\xi}\land y_{\zeta}\le x_{F({\xi},{\zeta})}.
\end{multline}
Then the topological density of $B$ is at least ${\mu}$.
\end{obs}

\begin{proof}[Proof of the observation]
Assume on the contrary that $\den(B)<{\mu}$.
Let ${\nu}$ be the minimal cardinal such that 
$X$ can be covered by  ${\nu}$ many ultrafilters,
$\{U_{\zeta}:{\zeta}<{\nu}\}$. Clearly ${\omega}\le{\nu}\le \den(B)<{\mu}$.

Let  
$F:\br {\nu}^+;2;\to {\nu^+}
$ be an injective function such that
 $(F({\alpha},{\gamma})\mod \nu)\ne (
 F({\beta},{\gamma}) \mod  \nu)$
for ${\alpha}<{\beta}<{\gamma}<{\nu}^+$. 
Then 
for each $I\in \br {\nu}^+;{\nu}^+;$ and ${\sigma}<{\nu}$
there
is $\{{\alpha},{\beta}\}\in \br I;2;$ such that 
${\sigma}< (F({\alpha},{\beta})\mod {\nu})$.

For each ${\xi}<{\nu}^+$  pick an element 
$x_{\xi}\in X\setm \bigcup\{U_{\zeta}:{\zeta}< ({\xi} \mod {\nu})\}$.
By (\ref{obs}) we can find 
$\<y_{\xi}:{\xi}<{\nu}^+\>\subs X$ such that 
$y_{\xi}\land y_{\zeta}\le x_{F({\xi},{\zeta})}$ for each
$\{{\xi},{\zeta}\}\in \br {\nu}^+;2;$.

Since $X\subs \bigcup_{{\zeta}<{\nu}}U_{\zeta}$, there is
${\rho}<{\nu}$ such that $I=\{{\xi}:y_{\xi}\in U_{\rho}\}$ 
is of size ${\nu}^+$.
Then there is $\{{\alpha},{\beta}\}\in \br I;2;$ such that 
${\rho}<(\ F({\alpha},{\beta}) \mod {\rho})$.
Thus $x_{F({\alpha},{\beta})}\notin U_{\rho}$ by the choice of
$x_{F({\alpha},{\beta})}$. On the other hand, 
$y_{\alpha}\land y_{\beta}\in U_{\rho}$, which contradicts
$y_{\alpha}\land y_{\beta} \le x_{F({\alpha},{\beta})}$.
\end{proof}

Unfortunately (\ref{obs}) is still 
too strong to be held in a c.c.c
Boolean algebra of topological density $\le {\mu}$. 
But, as it turns out,  it is not necessary to
consider all the sequences $\<x_{\xi}:{\xi}<{\nu}^+\>$
in (\ref{obs}) to derive 
that   the topological density of $B$ is at least ${\mu}$.
We will  introduce the index set $T$ in the proof below
in order to construct a manageable, but still 
large enough family of sequences.

\begin{proof}[A simplified proof of theorem \ref{tm:rs}]
The length of a sequence  ${\tau}$ is  denoted by $\len({\tau})$.
If $\len({\tau})={\alpha}+1$ then put $\base({\tau})={\tau}\restr {\alpha}$.
Given sequences ${\rho}$ and ${\tau}$ we write ${\rho}\ise{\tau}$
to mean that ${\rho}$ is an initial segment of ${\tau}$. 
Denote by $\trivalg$ the trivial  Boolean algebra $\{0,1\}$.

For each ${\nu}<{\mu}$ choose a function 
$h_{{\nu}^+}:\br {\nu}^+;2;\to {\nu^+}
$ such that
$h_{{\nu}^+}$ is 1--1 and  
$(h_{{\nu}^+}({\alpha},{\gamma})\mod \nu)\ne (
 h_{{\nu}^+}({\beta},{\gamma}) \mod  \nu)$
for ${\alpha}<{\beta}<{\gamma}<{\nu}^+$. 
Then 
for each $I\in \br {\nu}^+;{\nu}^+;$ and ${\sigma}<{\nu}$
there
is $\{{\alpha},{\beta}\}\in \br I;2;$ such that 
${\sigma}< (h_{{\nu}^+}({\alpha},{\beta})\mod {\nu})$.

\begin{definition}\label{d:T}
We define, by induction on ${\alpha}\le {\mu}^+$, $T_{\alpha}$ as follows:
\begin{enumerate}[(1)]%
\item $T_0=\{\emptyset\}$,
\item if ${\alpha}$ is limit, then 
$T_{\alpha}=\bigcup\{T_{\beta}:{\beta}<{\alpha}\}$,
\item if ${\alpha}={\beta}+1$ then let\newline  
\begin{tabular}{l}
$T_{\alpha}$= $T_{\beta}\cup\bigl\{{\rho}$:
${\rho} 
\text{ is a sequence of length $<{\mu}^+$
and for each ${\zeta}<\len({\rho})$}$\makebox[10mm]{}\\
\hfill either $\rho(\zeta)\in \mu$ or $\bigl( 
{\rho}({\zeta})\in T_{\beta}$ and 
$ {\rho}\restr {\zeta}\ise {\rho}({\zeta})\bigr)
\bigr\}$.
\end{tabular}
\end{enumerate}
\end{definition}

Let $T=T_{{\mu}^+}$.
For ${\eta}\in T$  let 
$\rank({\eta})=\min\{{\alpha}:{\eta}\in T_{\alpha}\}$.
For ${\rho}\in T$ we say that $\rho$ is a {\em \nt\ } iff
we can write ${\rho}=\hes{\rho};\, \conc\, \tas{\rho};$, where 
$(\len(\hes {\rho};)\mod {\mu})= 0 $, 
$\len(\tas\rho;)=\nu^+$ for some $\nu<\mu$ and 
$\rho\rest\zeta\ise \rho(\zeta)$ 
for each $\len(\hes \rho;)\le \zeta<\len(\rho)$.

Let $\tnt=\{\rho\in T: \text{$\rho$ is \nt}\}$. 
For  ${\rho}\in \tnt$
put  $E({\rho})=
\{{\rho}\conc \<i\>:i<\len({\rho})\}$.
The sets $E({\rho})$ are pairwise disjoint.
\begin{definition}
Define the function $F$ as follows. Let 
\begin{displaymath}
\dom(F)=\bigcup\{\br E({\rho});2;:{\rho}\in \tnt\}.
\end{displaymath}
For  $\{{\tau}_0,{\tau}_1\}\in \br E({\rho});2;$ write 
${\tau}_i={\rho}\conc \<k_i\>$ and  ${\rho}=\hes{\rho};\conc\tas{\rho};$
and put 
\begin{displaymath}
F({\tau}_0,{\tau}_1)=
\tas {\rho};(h_{\len({\tas\rho;})}(k_0,k_1))
\end{displaymath}
\end{definition}

\begin{definition}
For ${\eta}\in T$ let
\begin{displaymath}
D({\eta})=\bigl\{\ \{{\eta}\conc \<\emptyset\>\conc \<{\omega}n+i\>,
{\eta}\conc \<\emptyset\>\conc \<{\omega}n+j\>\}\ :\ i<j<n<{\omega}\bigr\}
\end{displaymath}
and $D=\bigcup\{D({\eta}):{\eta}\in T\}$.
\end{definition}

\begin{definition} 
Let  $\calbt$ be the Boolean algebra generated 
by $\{x_{\eta}:\eta\in T\}$ freely, except the relations in the 
following set  $\Gamma\cup\Delta$:
\begin{displaymath}
\Gamma=\bigl\{ 
x_{\tau_{0}}\land x_{\tau_{1}} \le x_{F(\tau_{0},\tau_{1})}
:\{{\tau}_0,{\tau}_1\}\in\dom(F) \bigr\},
\end{displaymath}
\begin{displaymath}
\Delta=\bigl\{x_{\tau_{0}}\land x_{\tau_{1}} =0
:\{{\tau}_0,{\tau}_1\}\in D \bigr\}.
\end{displaymath}
\end{definition}

\begin{lemma}\label{lm:large}
$\den(\calbt)\ge {\mu}$.
\end{lemma}

\begin{proof}
Let 
\begin{multline}\notag
{\nu}=\min\bigl\{{\nu}':\text{there are }  {\eta}\in T 
\text{ and ultrafilters }\<U_j:j<{\nu}'\>\\
\text{such that }\{x_{\tau}:{\eta}\ise {\tau}\}\subset
\bigcup\{U_j:j<{\nu}'\}\bigr\}.
\end{multline}
It is enough to show that ${\nu}={\mu}$. 
Assume on the contrary that
${\nu}<{\mu}$ witnessed by ${\eta}\in T$ and ultrafilters
$\<U_j:j<{\nu}\>$. We can assume that $\len({\eta})\equiv 0 \mod {\mu}$.
First observe that ${\nu}\ge {\omega}$ because  $D({\eta})$
can not be covered by finitely many ultrafilters either.

Construct a sequence
${\rho}$ of length ${\nu}^+$ such that for each ${\zeta}<{\nu}^+$:
\begin{equation}\tag{$\star$}\label{star}
\text{
${\eta}\conc ({\rho}\rest{\zeta})\ise {\rho}({\zeta})\in T$ and
$x_{{\rho}({\zeta})}\notin \bigcup \{U_j:j<({\zeta} \mod {\nu})\}$.
}
\end{equation}
Assume we have constructed ${\rho}({\xi})$ for ${\xi}<{\zeta}$
satisfying (\ref{star}). Then 
${\eta}\conc ({\rho}\restr {\zeta})\in T$ by definition.
Since ${\nu}$ was minimal
we have $$\{x_{{\tau}}:{\eta}\conc
({\rho}\rest{\zeta})\ise{\tau}\}\not\subset
\bigcup\bigl\{U_j:j<({\zeta}\mod{\nu})\bigr\}$$ and so we can choose a suitable
${\rho}({\zeta})$.

Now ${\eta}\conc {\rho}\in \tnt$.  For $i<{\nu}^+$
write ${\tau}_i={\eta}\conc {\rho}\conc \<i\>$. 
Since $\{{\tau}_i:i<{\nu}^+\}\subs \bigcup\{U_{\alpha}:{\alpha}<{\nu}\}$
there are $I\in \br {\nu}^+; {\nu}^+;$ and ${\alpha}<{\nu}$ such that 
${\tau}_i\in U_{\alpha}$ for each $i\in I$. Pick $\{i,j\}\in \br I; 2;$
such that $h_{{\nu}^+}(i,j)> ({\alpha} \mod {\nu})$. Then
$x_{F({\tau}_i,{\tau}_j)}=x_{{\rho}(h_{{\nu}^+}(i,j))}\notin
U_{\alpha}$ by (\ref{star}). 
But  $x_{{\tau}_i}\land x_{{\tau}_j}
\in U_{\alpha}$ which contradicts 
$x_{{\tau}_i}\land x_{{\tau}_j}\le x_{{F({\tau}_i,{\tau}_j)}}$.
\end{proof}

\begin{definition}\label{d:closed}
We say that $X\subs T$ is {\em closed}
provided:
\begin{enumerate}[(i)]%
\item \label{d:ise} if ${\rho}\in X$ and  $\len({\rho})={\alpha}+1$ then 
$\base({\rho})={\rho}\restr {\alpha}\in X$,
\item \label{d:F} if $\{{\tau}_1,{\tau}_2\}\in \dom (F)\cap \br X;2;$ then 
$F({\tau}_1,{\tau}_2)\in X$,
\item \label{d:D} if $\{{\tau}_1,{\tau}_2\}\in D$ and ${\tau}_1\in X$ then
${\tau}_2\in X$.
\end{enumerate}
\end{definition}

\begin{lemma}
Every  finite $X\subset T$ is contained in a  
finite closed $Y\subset T$.
\end{lemma} 
\begin{proof}
Observing $D\cap \dom(F)=\empt$ close $X$ first for $\ref{d:ise}$ and
$\ref{d:F}$, then close for 
$\ref{d:D}$. 
\end{proof}

\begin{definition}
Let $X\subset T$ be closed and $f:X\to\trivalg$. We say that 
$(*)_f$ holds iff 
\begin{enumerate}[(1)]
\item $f({\tau}_1)\land  f({\tau}_2)=0$ for each 
$\{{\tau}_1,{\tau}_2\}\in D\cap \br \dom(f);2;$,
\item  $f({\tau}_1)\land  f({\tau}_2)\le f(F({\tau}_1,{\tau}_2))$
for each $\{{\tau}_1,{\tau}_2\}\in \dom(F)\cap \br \dom(f);2;$.
\end{enumerate}
\end{definition}
The following lemma is a special case of  a well-known fact.
\begin{lemma}
Assume that  $f:T\to \trivalg$. Then there is a (unique) homomorphism 
${\varphi}_f$ from $\calbt$ into $\trivalg$ such that
${\varphi}_f(x_{\tau})=f({\tau})$  iff  $(*)_f$ holds.
\end{lemma}

For each $a\in \bcal_T\setm \{0\}$ fix a homomorphism
${\psi}_a:\bcal_T\to\trivalg$ with ${\psi}_a(a)=1$,
and finite, closed set $X_a\subs T$ such that 
$a$ is a Boolean combination of $\{x_{\tau}:{\tau}\in X_a\}$.
Define $f_a:X_a\to \trivalg$ by $f_a({\tau})=\psi_a(x_{\tau})$.

\begin{lemma}\label{lm:pre} $\calbt$  has precaliber ${\kappa}$ for each
$\kappa=\cf(\kappa) >\aleph_{0}$. 
\end{lemma}

\begin{proof}
Let $\{a_{\alpha}:{\alpha}<{\kappa}\}\subset \calbt\setminus\{0\}$.

It is enough to define a map $f:T\to \trivalg$ satisfying
$(*)_f$ such that $|\{{\alpha}:f_{a_{\alpha}}\subset f\}|={\kappa}$
because $f_{a_{\alpha}}\subset f$ implies
$1={\psi}_{a_{\alpha}}(a_{\alpha})={\varphi}_{f_{a_{\alpha}}}(a)=
{\varphi}_f(a_{\alpha})$. 

By thinning out $\{a_{\alpha}:{\alpha}<{\kappa}\}$ we can assume that
$\{X_{\alpha}: \alpha<\kappa \}$ is a $\Delta$-system with
kernel $X$ and that $f_{\alpha}\restr X=f'$.

\def\cp{crossing pair}
\def\cpf{\ccal\pcal}
A pair $\{{\tau}_1,{\tau}_2\}\in \dom(F)$ is called {\em \cp} if there are
${\alpha}\ne {\beta}<{\kappa}$ such that 
${\tau}_1\in X_{\alpha}\setminus X$ and ${\tau}_2\in
X_{\beta}\setminus X$. The family of crossing pairs is denoted my
$\cpf$.

\def\elso#1;{\bcal^{(1)}_{#1}}
\def\mas#1;{\bcal^{(2)}_{#1}}
\def\har#1;{\bcal^{(3)}_{#1}}
For ${\gamma}\in\oo$ let 
$$
\elso{\gamma};=\bigl\{\{{\tau}_o,{\tau}_1\}\in\cpf:
F({\tau}_0,{\tau}_1)\in X_{\gamma}\bigr\}.
$$
\begin{claim}\label{elso}
$|\elso{\gamma};|\le|X||X_{\gamma}|$.
\end{claim}

\begin{proof}[Proof of the claim \ref{elso}]
Let $\{{\tau}_0,{\tau}_1\}\in\elso {\gamma};$, 
${\tau}_i={\tau}\conc \<k_i\>$, 
${\tau}\in X$, 
${\eta}=F({\tau}_0,{\tau}_1)\in X_{\gamma}$.
Then the pair $\<{\tau},{\eta}\>$
determines the pair $\{{\tau}_0,{\tau}_1\}$.
Indeed, $F\restr \br E({\tau});2;$ is 1--1, so 
$\{{\tau}_0,{\tau}_1\}$ is the unique pair 
$\{{\tau}'_0,{\tau}'_1\}\in \br E({\tau});2;$ with 
$F({\tau}'_0, {\tau}'_1)={\eta}$.
\end{proof}

\def\tw{twins}
We say that ${\tau}_0$ and ${\tau}_1$ are {\em \tw} iff
${\tau}_0\ne{\tau}_1$ but $\base({\tau}_0)=\base({\tau}_1)$.
\begin{displaymath}
\mas{\gamma};=\bigl\{\{{\tau}_o,{\tau}_1\}\in\cpf:
\text{$\exists{\eta}\in X_{\gamma}$  
$F({\tau}_0,{\tau}_1)$ and ${\eta}$ are \tw }\bigr\}.
\end{displaymath}
\begin{claim}\label{mas}
$|\mas{\gamma};|\le|X||X_{\gamma}|$.
\end{claim}

\begin{proof}[Proof of claim \ref{mas}]
Let $\{{\tau}_0,{\tau}_1\}\in\mas {\gamma};$, 
${\tau}_i={\tau}\conc \<k_i\>$, ${\tau}\in X$. Fix ${\eta}\in X_{\gamma}$
such that $F({\tau}_0,{\tau}_1)$ and ${\eta}$ are twins.
Now the pair $\<{\tau},{\eta}\>$ determines 
the pair $\{{\tau}_0,{\tau}_1\}$.
Indeed, $F({\tau}_0,{\tau}_1)={\tau}({\xi})$ for some ${\xi}$, and there
is at most one ${\xi}$ such that ${\tau}({\xi})$ and ${\eta}$ are twins.
But ${\tau}$ and ${\xi}$ determine  $\{{\tau}_0,{\tau}_1\}$
because $F\restr \br E({\tau});2;$ is 1--1.
\end{proof}
Let
\begin{displaymath}
\har;=\bigl\{\ \bigl\{\{{\tau}_0,{\tau}_1\},\{{\tau}_2,{\tau}_3\}\bigr\}
\in\br \cpf;2;:
\text{ 
$F({\tau}_0,{\tau}_1)$ and $F({\tau}_2,{\tau}_3)$ are \tw }\bigr\}.
\end{displaymath}
\begin{claim}\label{har}
$|\har;|\le|X||X|$.
\end{claim}
\begin{proof}[Proof of claim \ref{har}]
Assume that 
$\bigl\{\{{\tau}_0,{\tau}_1\},\{{\tau}_2,{\tau}_3\}\bigr\} \in\br \cpf;2;$.
Then $\base({\tau}_0)=\base({\tau}_1)={\eta}\in X$ 
and $\base({\tau}_2)=\base({\tau}_3)={\rho}\in X$.
Moreover $F({\tau}_0,{\tau}_1)={\eta}({\xi})$ and 
$F({\tau}_2,{\tau}_3)={\rho}({\zeta})$ for some ${\xi}$ and ${\zeta}$.
But for given ${\eta}, {\rho}\in X$ there is at most one pair
$\{{\xi},{\zeta}\}$ such that ${\eta}({\xi})$ and ${\rho}({\zeta})$
are \tw. Since there is at most one 
$\{{\tau}'_0,{\tau}'_1\}\in E({\eta})$ with 
$F({\tau}'_0,{\tau}'_1)={\eta}({\xi})$ and
there is at most one 
$\{{\tau}'_2,{\tau}'_3\}\in E({\rho})$ with 
$F({\tau}'_2,{\tau}'_3)={\rho}({\xi}')$, we are done.
\end{proof}

So applying L\'az\'ar's free set mapping theorem 
we can thin out our sequence such that 
$\elso {\gamma};=\mas {\gamma};=\har;=\empt$ for each ${\gamma}\in {\kappa}$.

Let $f^-=\bigcup\{f_{\alpha}:{\alpha}<{\kappa}\}$.
Define $f:T\to\trivalg$ as follows.
Let $f({\eta})=1$ iff either  $f^-({\eta})=1$ or
 ${\eta}=F({\tau}_1,{\tau}_2\})$ for some $\{{\tau}_0,{\tau}_1\}\in \cpf$  
such that 
$f^-({\tau}_1)=f^-({\tau}_2)=1$.

Since $\bigcup\limits_{{\gamma}<{\kappa}}\elso {\gamma};=\empt$ 
we have $f^-\subset f$.

 We show that $(*)_f$ holds.
Assume that ${\tau}_0$ and ${\tau}_1$ are twins
and $f({\tau}_0)=f({\tau}_1)=1$.  
 Since $\har ;=\empt$ and
$\mas {\gamma};=\empt$ for each ${\gamma}<{\kappa}$, it follows 
that ${\tau}_0,{\tau}_1\in\bigcup\{X_{\alpha}:{\alpha}<{\kappa}\}=
\dom(f^-)$.
So $\{{\tau}_0,{\tau}_1\}\in D$ is impossible because 
${\tau}_0\in X_{\gamma}$ implies ${\tau}_1\in X_{\gamma}$.
Thus $\{{\tau}_0,{\tau}_1\}\in \dom (F)$ and 
$f^-({\tau}_0)=f^-({\tau}_1)=1$ and so $f(F({\tau}_0,{\tau}_1))=1$
by the construction of $f$.
\end{proof}

\begin{lemma}\label{lm:small}
$\den(\calbt)\le {\mu}$.
\end{lemma}

\begin{proof}
First fix a well-ordering $\prec$ of $T$ such that if
$\rank({\tau})<\rank({\tau}')$ then ${\tau}\prec {\tau}'$.

Let $T^-=T\setminus\bigcup\{E({\eta}):{\eta}\in \tnt\}$.
Consider the product space
\begin{displaymath}
\xcal=2^{T^-}\times (D_{[{\mu}]^{<{\omega}}})^{\tnt},
\end{displaymath}
where $D_{[{\mu}]^{<{\omega}}}$ denotes the discrete 
 topological space of size ${\mu}$ whose underlying set is 
$\br {\mu};<{\omega};$ instead of ${\mu}$.
Applying $\den ((D_{\mu})^{2^{\mu}})={\mu}$ and $|T|=2^{\mu}$ 
we can fix a dense
family $\{g_{\xi}:{\xi}<{\mu}\}\subs \xcal$. Write
$g_{\xi}=\<g^-_{\xi},g^*_{\xi}\>$.
For ${\xi}<{\mu}$ define $s_{\xi}:T\to \trivalg$ as follows:
$s_{\xi}\restr T^-=g^-_{\xi}$ and 
if ${\tau}\in T\setm T^-$, then pick the unique ${\eta}\in \tnt$
with ${\tau}\in E({\eta})$, ${\tau}={\eta}\conc \<i\>$, and let
$s_{\xi}({\tau})=1$ iff $i\in g^*_{\xi}({\eta})$.
Let 
\begin{displaymath}
\scal=\{s_{\xi}:{\xi}<{\mu}\}.	
\end{displaymath}
For ${\xi}<{\mu}$
define   $s^*_{\xi}:T\to \trivalg$  
by recursion  on $\prec$ as follows. 
Let $s^*_{\xi}({\tau})=1$ iff $s_{\xi}({\tau})=1$ and 
for each  $\{{\tau},{\tau}'\}\in D$ with  ${\tau}'\prec {\tau}$ 
we have  $s^*_{\xi}({\tau}')=0$ and
for each  $\{{\tau}',{\tau}\}\in \dom(F)$ with  ${\tau}'\prec {\tau}$ 
we have   $s^*_{\xi}({\tau}')\le s^*_{\xi}(F({\tau},{\tau}'))$. 

By induction on $\prec$ it is clear that  $(*)_{s^*_{\xi}}$ holds.

Now let $a\in \bcal_T\setm \{0\}$. By construction of 
$\scal$ we can find $s_{\xi}\in\scal$ such that 
$f_a\subs s_{\xi}$, moreover for each ${\eta}\in \tnt\cap X_a$
if $ {\tau}\in E({\eta})\setminus X_a$ then $s_{\xi}({\tau})=0$.

\begin{claim}
$s^*_{\xi}\supset f_a$.
\end{claim}
\begin{proof}
By induction on $\prec$.
Assume that the claim holds for ${\tau}'\in X_a$ provided 
${\tau}'\prec {\tau}$. We can assume $f_a({\tau})=1$.
If $\{{\tau}',{\tau}\}\in D$, ${\tau}'\prec{\tau}$ then ${\tau}'\in X_a$ so 
$f_a({\tau}')=0$ as $(*)_{f_a}$ holds. So, by the induction hypothesis,  
$s^*_{\xi}({\tau}')=f_a({\tau}')=0$.
Assume $\{{\tau}',{\tau}\}\in\dom (F)$, ${\tau}'\prec {\tau}$.
If ${\tau}'\in X_a$ then $F({\tau}',{\tau})\in X_a$ and 
$F({\tau}',{\tau})\prec {\tau}$ so $f_a({\tau}')=s^*_{\xi}({\tau}')$ and
$f_a(F({\tau}',{\tau}))=s^*_{\xi}(F({\tau}',{\tau}))$. 
Thus $s^*_{\xi}({\tau}')\le s^*_{\xi}(F({\tau}',{\tau}))=1$ because
$f_a({\tau}')=f_a({\tau})\land  f_a({\tau}')\le f_a(F({\tau}',{\tau}))$ as
$(*)_{f_a}$ holds. If ${\tau}'\notin X_a$ then $s_{\xi}({\tau}')=0$ 
by the assumption 
about $s_{\xi}$ and so $s^*_{\xi}({\tau}')=0$.
Thus $s^*_{\xi}({\tau})=1$ by the construction of $f^*$.
\end{proof}
Thus $f_a(a)=s^*_{\xi}(a)$, i.e.
$\bcal_T\setm
\{\empt\}=\bigcup\{{\varphi}_{s^*_{\xi}}^{-1}\{1\}:{\gamma}<{\mu}\}$,
which was to be proved, so the lemma holds.
\end{proof}
$\calbt$  is c.c.c by lemma \ref{lm:pre} and 
$\den(\calbt)={\mu}$ by lemmas  \ref{lm:small} and \ref{lm:large}
so the theorem is proved.
\end{proof}

\end{document}